\documentclass[lettersize,journal]{IEEEtran}
\usepackage{amsmath,amsfonts, units}
\usepackage{algorithmic}
\usepackage{algorithm}
\usepackage{array}
\usepackage[caption=false,font=normalsize,labelfont=sf,textfont=sf]{subfig}
\usepackage{textcomp}
\usepackage{stfloats}
\usepackage{url}
\usepackage{verbatim}
\usepackage{graphicx}
\usepackage{cite}
\usepackage{amsthm}
\usepackage{amssymb}
\usepackage{subfig}
\usepackage{booktabs, multirow}
\usepackage{xcolor}
\DeclareMathOperator*{\argmin}{arg\,min}

\newtheorem{theorem}{Theorem}[section]
\newtheorem{lemma}[theorem]{Lemma}
\newtheorem{assumption}{Assumption}
\theoremstyle{definition}
\newtheorem{example}{Example}
\theoremstyle{remark}
\newtheorem{remark}{Remark}
\begin{document}

\title{Stochastic Sequential Feedback Optimization for Wind Farm Control under Time-Varying Inflow}

\author{Shijie~Huang~and~Sergio~Grammatico
\thanks{S. Huang and S. Grammatico are with the Delft Center of Systems and Control, TU Delft, The Netherlands. E-mail address: \texttt{S.Huang-5@tudelft.nl; s.grammatico@tudelft.nl}.}
\thanks{This work was partially supported by the NWO under research project Online Optimization for Offshore Wind Farms.}}



\maketitle

\begin{abstract}
We consider steady-state optimization of nonlinear dynamical systems under time-varying stochastic operating conditions, where the objective is an expected-value function over the current operating-condition distribution, so that changes in the distribution induce a moving stochastic optimum. To address this problem, we propose stochastic sequential feedback optimization (S2FO), which combines projected feedback-gradient updates with sequential sensitivity estimates obtained by linearizing the plant at realized stochastic operating conditions. We analyze the closed-loop performance of S2FO through expected dynamic regret, instantaneous tracking in expectation, and almost-sure convergence under diminishing step sizes. We then apply the method to wind farm power maximization under time-varying stochastic inflow. Simulations on a benchmark layout and on a realistic offshore wind farm show that S2FO improves average power generation when persistent wake interactions are present, and that its advantage decreases as wind-direction variability increases.
\end{abstract}

\begin{IEEEkeywords}
Steady-state optimization, time-varying, stochastic feedback optimization, wind farm.
\end{IEEEkeywords}

\section{Introduction}
Many engineering systems are operated by repeatedly adjusting control inputs to improve their steady-state performance. This problem is known as steady-state optimization, where the goal is to find an operating point that optimizes a performance criterion of the controlled plant \cite{hauswirth2024optimization}. In many applications, the steady-state performance depends on external operating conditions, such as load variations in distribution grids~\cite{bernstein2019real}, process disturbances in chemical plants \cite{krishnamoorthy2022real}, or weather variability in renewable energy systems~\cite{franccois2014real}. When these operating conditions are random, the relevant objective is the expected steady-state cost over their distribution, rather than the cost evaluated at a nominal operating condition. If the underlying distribution also evolves over time, the stochastic optimum becomes a moving target that the optimizer must track in closed loop with the plant dynamics.

Wind farm control is a representative example of this setting, with the wind flow playing the role of the external operating condition \cite{meyers2022wind}. The goal is to improve farm-level power production by coordinating turbine inputs, such as yaw angles and axial induction factors. A standard greedy strategy maximizes the power of each turbine individually, but this can be suboptimal at the farm level because upstream turbines generate wakes that reduce the inflow velocity of downstream turbines~\cite{boersma2017tutorial}. Moreover, these wake interactions vary with wind speed, wind direction, and atmospheric stability~\cite{porte2020wind}. Because these inflow quantities are themselves time-varying and stochastic, the farm-level optimal operating point also changes over time~\cite{doekemeijer2020closed,howland2022optimal,becker2025closed}. Therefore, wind farm control provides a concrete motivation for studying how to optimize a dynamical system at steady state under time-varying stochastic operating conditions.

Existing wind farm control methods have addressed changing inflow in several ways. Robust wake-steering methods based on engineering wake models optimize yaw set-points over wind-direction variability or model uncertainty, so that the selected set-points are less sensitive to estimation and modeling errors~\cite{rott2018robust,simley2020design,quick2020wake}. Measurement-driven wake-steering methods use turbine measurements to refine inflow estimates or wake-model parameters and recompute yaw commands as conditions change~\cite{doekemeijer2020closed,howland2022optimal}. To account more explicitly for wake propagation, another line of work uses dynamic wake models within model predictive control. Adjoint-based MPC with WFSim has been tested under prescribed wind-direction changes~\cite{vali2019adjoint}. More recent work uses free-vortex wake models under designed and measured wind variations~\cite{van2024dynamic}, or combines FLORIDyn with flow-state estimation for MPC under time-varying inflow~\cite{becker2025closed}. These MPC methods explicitly account for wake transients over a finite prediction horizon. This paper adopts a complementary viewpoint based on stochastic steady-state optimization, where changing inflow is represented through a time-varying expected objective for a nonlinear plant.

Online feedback optimization (FO) is a natural framework for steady-state optimization in dynamical systems. It implements optimization updates directly in closed loop with the plant by using real-time output measurements and sensitivity information~\cite{hauswirth2024optimization}. Existing feedback optimization results have established stability and tracking guarantees when the steady-state sensitivity is available~\cite{colombino2019online,hauswirth2020timescale}. To reduce the reliance on exact sensitivities, model-free variants estimate gradient information directly from data \cite{he2023model}, while gray-box variants further combine approximate model information with measurement-based updates~\cite{he2026gray}. Recent work has also studied time-varying environments from an output-regulation perspective, showing that exact tracking in deterministic settings requires an internal model of the temporal variability~\cite{bianchin2025feedback}. Thus, the existing FO methods mainly address fixed or deterministically varying operating conditions, but they do not explicitly target stochastic operating conditions whose distributions evolve over time. Online stochastic optimization provides a related viewpoint for problems with time-varying distributions~\cite{cao2020online,drusvyatskiy2023stochastic}. In this direction, Bianchin et al.~\cite{bianchin2023online} linked stochastic optimization with dynamical systems by studying unknown linear systems and proposing a data-driven controller based on estimating the steady-state transfer function. The nonlinear case introduces an additional difficulty. The steady-state input-output sensitivity changes with the operating point and the external condition, and therefore cannot be represented by a single fixed gain. This applies in particular to wind farm control, where power sensitivities depend on the current flow state and wind condition.

Motivated by this gap, we develop stochastic sequential feedback optimization (S2FO) for nonlinear plants under time-varying stochastic operating conditions. The closest algorithmic precursor is deterministic SFO~\cite{huang2025sequential}, where the steady-state sensitivity is updated by repeatedly linearizing the nonlinear plant around the current operating point. The stochastic setting changes the optimization target itself. Instead of seeking the optimum associated with a fixed operating condition, the controller aims to track a moving stochastic optimum defined by an evolving distribution. The proposed S2FO evaluates the linearization at realized stochastic operating conditions and uses the resulting model-based sensitivity information in the feedback update. This distinction also changes the performance analysis. The closed-loop error now depends not only on the plant response and the sequential linearization error, but also on the randomness of the operating condition and the drift of its distribution.

The main contributions of this paper are as follows.

\begin{itemize}
  \item We propose a stochastic sequential feedback optimization algorithm for time-varying stochastic steady-state optimization of nonlinear dynamical systems. The algorithm combines projected feedback-gradient steps with sequential sensitivity updates evaluated at the realized condition of each iteration, thereby extending the deterministic paradigm of \cite{huang2025sequential} to a stochastic setting.
  \item We establish theoretical guarantees for S2FO, including an expected dynamic regret bound, an instantaneous tracking error bound in expectation, and almost-sure convergence under diminishing step sizes and asymptotically vanishing parameter variation. These results characterize how the closed-loop performance depends on the drift of the stochastic optimum and on the variation of the external operating condition.
  \item We apply S2FO to wind farm power maximization under time-varying stochastic wind conditions. Simulations on a benchmark \(2\times 3\) layout and on a realistic Lillgrund subset show that S2FO improves average farm power when persistent wake interactions are present. The results also show that S2FO maintains a positive gain over the greedy baseline across the tested statistical mismatch range, while its gain decreases when strong wind-direction variability raises the greedy baseline.
 \end{itemize}

\section{Problem Formulation and Algorithm}
\label{sec:problem}

We consider a discrete-time nonlinear dynamical system
\begin{equation}
\label{eq:dynamics}
\begin{cases}
x_{t+1} = f(x_t, u_t, \theta_t)  \\
y_t = g(x_t, u_t, \theta_t)
\end{cases}
\end{equation}
where $x_t \in \mathbb{R}^n$ is the system state, $u_t \in \mathcal{U} \subseteq \mathbb{R}^p$ is the control input, $y_t \in \mathbb{R}^m$ is the measured output, and $\theta_t \in \Theta \subseteq \mathbb{R}^d$ is a stochastic parameter. We assume that $\mathcal{U}$ is a compact and convex set that represents physical constraints on the control inputs.

In wind farm applications, $\theta_t$ captures time-varying atmospheric conditions such as wind speed, direction, and turbulence intensity, which evolve according to meteorological processes. The state $x_t$ represents the spatially discretized wind velocity field, $u_t$ represents control inputs including axial induction factors and yaw angles, and $y_t$ represents performance outputs such as power generation~\cite{boersma2017tutorial}.

\vspace{-2mm}
\subsection{Steady-State Mapping and Stochastic Optimization}

For each given parameter realization $\bar{\theta} \in \Theta$ and control input $\bar{u} \in \mathcal{U}$, we assume the system admits a unique steady state. Specifically, there exists a steady-state mapping $\phi: \mathcal{U} \times \Theta \to \mathbb{R}^n$ such that
\begin{equation}
\phi(\bar{u}, \bar\theta) = f(\phi(\bar u, \bar\theta), \bar u, \bar\theta), \quad \forall \bar u \in \mathcal{U}, \bar\theta \in \Theta.
\end{equation}
The steady-state input-output mapping is given by
\begin{equation}
\label{eq:ss_mapping}
\bar y = g(\phi(\bar u, \bar\theta), \bar u, \bar\theta) =: h(\bar u, \bar\theta).
\end{equation}

Let $(\theta_t)_{t \in\mathbb{N}}$ be a stochastic process taking values in $\Theta$, and let $P_t\in\mathcal P_1(\Theta)$ denote the marginal distribution of $\theta_t$, where $\mathcal P_1(\Theta)$ is the set of probability measure on $\Theta$ with finite first moment. For each time $t$, the stochastic steady-state cost associated with an input $u\in\mathcal{U}$ is given by
\begin{equation}
\label{eq:F_t}
F_t(u) := \mathbb{E}_{\theta_t \sim P_t}[J(u, h(u, \theta_t))],
\end{equation}
where $J: \mathcal{U} \times \mathbb{R}^m \to \mathbb{R}$ is the cost function. The corresponding stochastic steady-state optimization problem is
\begin{equation}
\label{eq:stoch_opt}
\min_{\bar{u}\in\mathcal{U}} \quad F_t(\bar{u}).
\end{equation}
Let $u_t^\star$ denote an optimal solution of \eqref{eq:stoch_opt}. Since $P_t$ in \eqref{eq:F_t} may vary over time, the optimizer $u_t^\star$ is generally time-varying. The goal is therefore to generate a feasible input sequence $(u_t)_{t \in\mathbb{N}}$ that tracks $(u_t^\star)_{t \in\mathbb{N}}$ in closed loop.

We make the following assumptions on the system dynamics and cost function, which extend Assumptions~1 and~2 in \cite{huang2025sequential} to hold uniformly over the parameter space $\Theta$.

\begin{assumption}
\label{ass:system}
The system dynamics in \eqref{eq:dynamics} satisfy the following conditions uniformly over $\theta \in \Theta$:
\begin{enumerate}
\item The mapping $x \mapsto f(x, u, \theta)$ is uniformly contractive with factor $\rho_f \in (0, 1)$:
\begin{equation}
\sup_{(x,u,\theta) \in \mathbb{R}^n \times \mathcal{U} \times \Theta} \|\nabla_x f(x, u, \theta)\| \leq \rho_f < 1.
\end{equation}

\item The partial gradients $\nabla_x f$, $\nabla_u f$, and $\nabla_\theta f$ are uniformly bounded and Lipschitz continuous.

\item The steady-state mapping $h: \mathcal{U} \times \Theta \to \mathbb{R}^m$ in \eqref{eq:ss_mapping} is Lipschitz continuous in both arguments, i.e., 
\begin{align*}
\|h(u_1, \theta) - h(u_2, \theta)\| &\leq L_{h,u} \|u_1 - u_2\|, \\
\|h(u, \theta_1) - h(u, \theta_2)\| &\leq L_{h,\theta} \|\theta_1 - \theta_2\|,
\end{align*}
for some constants $L_{h,u}$ and $L_{h,\theta}$.
\end{enumerate}
\end{assumption}

\begin{assumption}
\label{ass:cost}
The cost function $J: \mathcal{U} \times \mathbb{R}^m \to \mathbb{R}$ in \eqref{eq:F_t} satisfies:
\begin{enumerate}
\item The partial gradients $\nabla_u J$ and $\nabla_y J$ are uniformly bounded and Lipschitz continuous with constants $G_u^J$, $G_y^J$, $L_{J,u}$, and $L_{J,y}$.
\item For any fixed $y \in \mathbb{R}^m$, the mapping $u \mapsto J(u, y)$ is $\mu_J$-strongly convex:
\begin{equation*}
(\nabla_u J(u_1, y) - \nabla_u J(u_2, y))^\top (u_1 - u_2) \geq \mu_J \|u_1 - u_2\|^2.
\end{equation*}
\end{enumerate}
\end{assumption}

The uniformity in $\theta$ is needed for the subsequent convergence and tracking analysis, where the same Lipschitz and contractivity constants are used across all realizations of $\theta_t$.

\subsection{Stochastic Sequential Feedback Optimization}

Standard feedback optimization algorithms for nonlinear plants require the steady-state input-output sensitivity $\nabla h$, which is generally not available in closed form for high-dimensional plants such as wind farms. To address this difficulty, the deterministic sequential feedback optimization (SFO) algorithm of~\cite{huang2025sequential} updates the linearization point at every iteration. At each step, the system is linearized around the current operating point $(\hat x_t, \hat u_t)$, and the resulting linearized model yields an approximate sensitivity $\widetilde{\nabla h}(\hat x_t, \hat u_t)$. A projected feedback-gradient step is then taken using this approximate sensitivity together with the measured output. As the iterates approach the optimal steady state, the linearization error decreases, and the algorithm converges to a neighborhood of the optimum.

The stochastic time-varying setting in~\eqref{eq:stoch_opt} introduces a new challenge. The operating condition $\theta_t$ is no longer fixed or known in advance. Instead, it is modeled as a stochastic process whose distribution $P_t$ evolves over time. As a result, the optimum $u_t^\star$ becomes a moving target induced by evolving statistics, and the algorithm aims to track this target rather than converge to a fixed point. In this paper, we propose the stochastic sequential feedback optimization (S2FO) algorithm, given in Algorithm~\ref{alg:ssfo}, which extends the SFO framework to this setting by evaluating the sequential linearization at the realized $\theta_t$ of each iteration.

\begin{algorithm}[H]
\caption{Stochastic Sequential Feedback Optimization (S2FO)}
\label{alg:ssfo}
\begin{algorithmic}[1]
\REQUIRE Initial state $x_0$, initial input $u_0 \in \mathcal{U}$, step size $\alpha > 0$

\FOR{$t = 0, 1, 2, \ldots$}
\STATE \textbf{Sample}: $\theta_t\sim P_t$.
\STATE \textbf{Output measurement}: $y_t = g(x_t, u_t, \theta_t)$.
\STATE \textbf{Sensitivity estimation:} linearize \eqref{eq:dynamics} at $(x_t, u_t, \theta_t)$ and compute
\begin{equation}
\widetilde{\nabla h}(x_t, u_t, \theta_t) = \nabla_x g \cdot (I - \nabla_x f)^{-1} \nabla_u f + \nabla_u g
\end{equation}
with all gradients evaluated at $(x_t, u_t, \theta_t)$.
\STATE \textbf{Stochastic gradient computation:}
\begin{equation}\label{eq:gradient}
g_t = \nabla_u J(u_t, y_t) + \widetilde{\nabla h}(x_t, u_t, \theta_t)^\top \nabla_y J(u_t, y_t)
\end{equation}
\STATE \textbf{System evolution:} $x_{t+1} = f(x_t, u_t, \theta_t)$. 
\STATE \textbf{Projected gradient update:}
\begin{equation}\label{eq:update}
u_{t+1} = \Pi_{\mathcal{U}}\left(u_t - \alpha g_t\right)
\end{equation}
\ENDFOR
\end{algorithmic}
\end{algorithm}

\begin{remark}
\label{rem:mean_based_update}
A natural deterministic alternative to Algorithm~\ref{alg:ssfo} is to replace the random operating condition by its mean $\bar\theta_t := \mathbb{E}[\theta_t]$ and compute the sensitivity at $\bar\theta_t$. This mean-based update is analogous to a certainty-equivalence approximation, where deterministic estimates are used in place of stochastic quantities in the controller design~\cite{van1981certainty}. It recovers the deterministic SFO update in~\cite{huang2025sequential} applied at the mean trajectory $\bar\theta_t$. However, this update is generally not equivalent to optimizing the stochastic steady-state objective. Defining $\bar F_t(u) := J(u, h(u, \bar\theta_t))$, we have $F_t(u) \neq \bar F_t(u)$ and $\nabla F_t(u) \neq \nabla \bar F_t(u)$ in general, because $h(u, \cdot)$ in \eqref{eq:ss_mapping} is nonlinear at the operating condition. Hence, mean-based sensitivity information may miss the effect of stochastic variations on the expected steady-state performance. This distinction motivates the comparison between stochastic and mean-based sensitivity updates in the numerical study.
\end{remark}

Let $\mathcal{F}_t := \sigma(\theta_1, \ldots, \theta_t, u_1, \ldots, u_t, y_1, \ldots, y_t)$ denote the $\sigma$-algebra generated by all observations up to time $t$. We make the following assumption on the stochastic parameter.

\begin{assumption}
\label{ass:cond_indep}
For each \(t\in\mathbb{N}\), the input \(u_t\) is \(\mathcal{F}_{t-1}\)-measurable, and the stochastic parameter \(\theta_t\) is drawn independently of \(\mathcal{F}_{t-1}\) from the current distribution \(P_t\) that defines the objective in \eqref{eq:F_t}. 
\end{assumption}

Assumption~\ref{ass:cond_indep} is a standard sampling assumption in stochastic approximation and online stochastic optimization~\cite{nemirovski2009robust,cao2020online}. It says that the control input in \eqref{eq:dynamics} is chosen using past information before $\theta_t$ is sampled from \(P_t\). This assumption is used to establish conditional unbiasedness of the realized steady-state gradient, which is needed in the expected dynamic-regret analysis.

\subsection{Comparison with Closely Related Algorithms}\label{sec:comparison}

Let us compare S2FO with three closely related lines of work to clarify the formulation and algorithmic differences.

\paragraph*{Online stochastic optimization for unknown linear systems} Bianchin et al.~\cite{bianchin2023online} consider a problem setup that is structurally close to~\eqref{eq:stoch_opt}, namely the regulation of a stochastic linear dynamical system to the solution of a stochastic convex optimization problem under time-varying disturbances. In that setting, the plant is linear, so the steady-state input-output map is captured by a single matrix gain that can be estimated from input-output trajectory data. S2FO instead addresses nonlinear plants whose steady-state sensitivity changes with the operating point and the parameter realization, so a single fixed gain is no longer sufficient. This motivates the sequential sensitivity update in Algorithm~\ref{alg:ssfo}.

\paragraph*{Gray-box feedback optimization for nonlinear systems} He et al.~\cite{he2026gray} address nonlinear plants by combining model-based sensitivities with model-free updates through an adaptive convex combination. They characterize the closed-loop performance through dynamic regret in a constrained time-varying setting, but the disturbance sequence in their formulation is deterministic. In contrast, our S2FO faces a stochastic $\theta_t$ with a time-varying distribution. This stochastic structure motivates the use of the realized $\theta_t$ in the sensitivity update (Remark~\ref{rem:mean_based_update}) and complicates the associated closed-loop performance analysis.

\paragraph*{Online stochastic optimization with time-varying distributions} Cao et al.~\cite{cao2020online} study expected objective functions that change over time and provide dynamic-regret guarantees, but do not include a dynamical plant. In S2FO, the stochastic gradient information is constructed from the plant output and a model-based steady-state sensitivity, which introduces additional output mismatch and sensitivity approximation errors.

\section{Main Results}
In this section, we establish theoretical guarantees for the proposed stochastic sequential feedback optimization algorithm. We present three main results: (i) an expected dynamic regret bound, (ii) a stochastic instantaneous tracking error bound, and (iii) an almost-sure convergence under a weighted summability condition on the parameter drift. To simplify the notation, we assume throughout this section that $g(x, u, \theta) = x$, so that the output coincides with the state. The analysis extends to general output functions $g$ with Lipschitz gradients, introducing only additional notational complexity.

Let us first establish preliminary lemmas that characterize the sensitivity approximation, output mismatch, and gradient estimation error.
\subsection{Preliminary Results}\label{sec:pre_lemmas}
\begin{lemma}
\label{lem:sensitivity_lip}
Suppose Assumption~\ref{ass:system} holds. The approximate sensitivity map $(x, u, \theta) \mapsto \widetilde{\nabla h}(x, u, \theta)$ defined by $\widetilde{\nabla h}(x, u, \theta) = (I - \nabla_x f(x, u, \theta))^{-1} \nabla_u f(x, u, \theta)$ is Lipschitz continuous. Specifically, for all $(x_1, u_1, \theta_1), (x_2, u_2, \theta_2)$:
\begin{multline*}
\|\widetilde{\nabla h}(x_1, u_1, \theta_1) - \widetilde{\nabla h}(x_2, u_2, \theta_2)\|\\ \leq C (\|x_1 - x_2\| + \|u_1 - u_2\| + \|\theta_1 - \theta_2\|),
\end{multline*}
where $C = \frac{(1-\rho_f)L_{f,u} + G_u^f L_{f,x}}{(1-\rho_f)^2}$ and $G_u^f := \sup_{x,u,\theta} \|\nabla_u f(x,u,\theta)\|$.
\end{lemma}

In Lemma~\ref{lem:sensitivity_lip}, the constant \(C\) increases as \(\rho_f\) approaches one. This means that when the plant contracts weakly, even small perturbations in the state or operating condition can cause large changes in the sensitivity approximation.

We next analyze how the plant output $y_t$ deviates from the steady-state output $h(u_t,\theta_t)$ associated with the current input and operating condition. 
\begin{lemma}
\label{lem:output_error}
Under Assumption~\ref{ass:system}, there exist constants $C_u, C_\theta > 0$ such that the output mismatch
\begin{equation}\label{eq:output_mismatch}
e_t := \|y_t - h(u_t, \theta_t)\|
\end{equation}
satisfies
\begin{equation}
\label{eq:e_recursion}
e_{t+1} \leq \rho_f e_t + C_u \|u_{t+1} - u_t\| + C_\theta \|\theta_{t+1} - \theta_t\|.
\end{equation}
Moreover, for any $T \geq 1$,
\begin{align*}
\sum_{t=1}^{T} e_t &\leq \frac{e_1}{1 - \rho_f} + \frac{C_u}{1 - \rho_f} \sum_{t=1}^{T-1} \|u_{t+1} - u_t\|\\
&\quad + \frac{C_\theta}{1 - \rho_f} \sum_{t=1}^{T-1} \|\theta_{t+1} - \theta_t\|.
\end{align*}
\end{lemma}

The bound in Lemma~\ref{lem:output_error} consists of three terms. The first term depends on the initial output mismatch $e_1$ and is a constant that does not grow with $T$. The other two terms reflect the variation of the control input and of the stochastic operating condition over the horizon.

We now use the output mismatch in \eqref{eq:output_mismatch} to bound the gradient estimation error.

\begin{lemma}
\label{lem:grad_error}
Let $g_t$ denote the gradient estimate in~\eqref{eq:gradient} and let $\tilde{J}(u, \theta) := J(u, h(u, \theta))$, with $h$ as in \eqref{eq:ss_mapping}, denote the steady-state cost for a given $\theta$. Under Assumptions~\ref{ass:system}--\ref{ass:cost}, there exists a constant $\bar L > 0$ such that
\begin{equation}
  \label{eq:grad_mismatch_bound}
  \|\nabla \tilde J(u_t,\theta_t)-g_t\|\le \bar L e_t.
  \end{equation}
  Moreover, under Assumption~\ref{ass:cond_indep},
  \begin{equation}
  \label{eq:cond_unbiased}
  \mathbb{E}\!\left[\nabla \tilde J(u_t,\theta_t)\mid \mathcal F_{t-1}\right]
  =
  \nabla F_t(u_t).
  \end{equation}
\end{lemma}

Lemma~\ref{lem:grad_error} relates the gradient estimate \(g_t\) to the gradient \(\nabla F_t(u_t)\) of the expected steady-state cost through the intermediate quantity \(\nabla \tilde J(u_t,\theta_t)\). The first inequality bounds the deviation between \(g_t\) and \(\nabla \tilde J(u_t,\theta_t)\), which arises because the plant has not yet reached the steady state. This deviation is controlled by the output mismatch $e_t$. The second equality states that \(\nabla \tilde J(u_t,\theta_t)\) is an unbiased estimate of \(\nabla F_t(u_t)\) given $\mathcal{F}_{t-1}$. We denote the corresponding zero-mean noise by $\zeta_t:=\nabla \tilde J(u_t,\theta_t) - \nabla F_t(u_t)$, which appears throughout the subsequent analysis.

\subsection{Expected Dynamic Regret}\label{sec:regret}

We now bound the expected dynamic regret of Algorithm~\ref{alg:ssfo}, which measures the cumulative suboptimality of the iterates $u_t$ relative to the time-varying optimum $u_t^{\star}$:
\begin{equation}
\mathrm{reg}(T) := \sum_{t=1}^{T} \big(F_t(u_t) - F_t(u_t^{\star})\big).
\end{equation}
To derive the bound, we adopt two standard non-stationarity measures from the dynamic regret literature~\cite{mokhtari2016online,cheung2019learning}: the path length of the optimum $\Delta_T := \sum_{k=1}^{T-1} \|u_{k+1}^\star - u_{k}^\star\|$, which captures the temporal variation of the moving target, and the variation budget of the operating condition $V_T := \sum_{k=1}^{T-1} \mathbb{E}[\|\theta_{k+1} - \theta_k\|]$, which captures the cumulative drift of the parameter in expectation. 
\begin{theorem}
\label{thm:dynamic_regret}
Suppose Assumptions~\ref{ass:system}--\ref{ass:cond_indep} hold and that the expected steady-state costs $F_t$ in \eqref{eq:F_t} are convex on $\mathcal U$. Let $R > 0$ be such that $\sup_{u\in\mathcal{U}}\|u\| \leq R$. For any $T \geq 1$ and constant step size $\alpha > 0$, the expected dynamic regret of Algorithm~\ref{alg:ssfo} satisfies
\begin{equation}\label{eq:regret_bound}
\begin{aligned}
\mathbb{E}[\mathrm{reg}(T)] &\leq \frac{1}{2\alpha}\big(4R^2 + 4R \Delta_T\big) + \frac{\alpha G^2}{2} T \\
&\quad + \frac{2R \bar{L}}{1 - \rho_f} \left( e_1 + C_u \alpha G T + C_\theta V_T \right),
\end{aligned}
\end{equation}
for some constant $G > 0$. 
\end{theorem}

The bound in Theorem~\ref{thm:dynamic_regret} consists of three parts. The first two terms are standard in dynamic regret analyses of projected gradient methods \cite[Theorem 3]{cheung2019learning}. The last term comes from the cumulative output mismatch $\sum_t e_t$ analyzed in Lemma~\ref{lem:output_error}.
\begin{remark}\label{rem:step_size}
When both the path length and the parameter variation budget grow sublinearly, i.e., $\Delta_T = o(T)$ and $V_T = o(T)$, choosing $\alpha = \sqrt{(1+\Delta_T) / T}$ yields
\begin{equation*}
\mathbb{E}[\mathrm{reg}(T)] = \mathcal{O}\big(\sqrt{T(1 + \Delta_T)}\big) + \mathcal{O}(V_T) + \mathcal{O}(1) = o(T).
\end{equation*}
Theorem~\ref{thm:dynamic_regret} shows that S2FO achieves sublinear expected dynamic regret when the stochastic optimum and the operating condition do not vary too rapidly.
\end{remark}

\begin{example}\label{ex:wind_decomp}
To illustrate a scenario where the parameter variation budget satisfies $V_T = o(T)$, consider an additive decomposition of the stochastic parameter, for example the wind condition~\cite{zarate2013continuous}:
\begin{equation}\label{eq:theta_decomp}
\theta_t = m_t + \sigma_t \xi_t,
\end{equation}
where $m_t \in \mathbb{R}^d$ represents a deterministic trend (e.g., daily wind patterns), $\xi_t \in \mathbb{R}^d$ is a zero-mean stochastic process with bounded expected norm $\mathbb{E}[\|\xi_t\|] \leq C_{\xi}$, and $\sigma_t \in \mathbb{R}_{\geq 0}$ is a time-varying noise magnitude. The expected one-step variation of $\theta_t$ satisfies
\begin{equation}\mathbb{E}[\|\theta_{t+1} - \theta_t\|] \leq \|m_{t+1} - m_t\| + (\sigma_{t+1} + \sigma_t)C_{\xi}.\end{equation}
A sufficient condition for $V_T  = o(T)$ is that:
\begin{enumerate}
\item The deterministic trend has sublinear cumulative drift: $\sum_{t=1}^T \|m_{t+1} - m_t\| = o(T)$.
\item The fluctuation magnitude decays as $\sigma_t = \mathcal{O}(t^{-\gamma})$ for some $\gamma \in (0, 1]$.
\end{enumerate}
\end{example}

\subsection{Stochastic Instantaneous Tracking Error}\label{sec:instan_tracking}

While dynamic regret provides a cumulative performance measure, instantaneous guarantees on the tracking error at each time step are also of interest. Unlike the regret analysis, where the stochastic gradient noise $\zeta_t$ vanishes in conditional expectation, the instantaneous bound must explicitly control its magnitude through a second-moment condition. We therefore impose the following assumption, analogous to the standard bounded-variance assumption in \cite[Assumption A1(b)]{ghadimi2013stochastic}.

\begin{assumption}
  \label{ass:noise}
  There exists a constant $\sigma_g>0$ such that, for all $t\geq 1$ and all $u\in\mathcal U$,
  \begin{equation}
  \mathbb E\!\left[
  \left\|
  \nabla\tilde J(u,\theta_t)-\nabla F_t(u)
  \right\|^2
  \mid \mathcal F_{t-1}
  \right]
  \leq \sigma_g^2.
  \end{equation}
  \end{assumption}

  Assumption~\ref{ass:noise} bounds the conditional second moment of the realized steady-state gradient noise. In particular, Jensen's inequality implies $\mathbb E[\|\zeta_t^\star\|]\leq \sigma_g$, where $\zeta_t^\star :=
  \nabla\tilde J(u_t^\star,\theta_t)-\nabla F_t(u_t^\star)$ is the noise evaluated at the optimum.

We now establish a bound on the instantaneous tracking error $e_t^u := \|u_t - u_t^{\star}\|$. Define the coefficient matrix
\begin{equation}
\label{eq:M_matrix}
M := \begin{bmatrix} \sqrt{1 - 2\alpha\mu_J + \alpha^2 L_{J,u}^2} + \alpha C_1 & \alpha C_2 \\ G_u^f & \rho_f \end{bmatrix}
\end{equation}
where $C_1 := L_{h,u} L_{J,y} + G_y^J C$, $C_2 := L_{J,u} + L_{h,u} L_{J,y}$ are constructed from the Lipschitz constants in Assumptions~\ref{ass:system} and~\ref{ass:cost}, and $C$ is defined in Lemma~\ref{lem:sensitivity_lip}.

\begin{theorem}
\label{thm:tracking_error}
Suppose Assumptions~\ref{ass:system}--\ref{ass:noise} hold, and that the step size $\alpha > 0$ is chosen such that $\rho(M) < 1$, with $M$ as defined in \eqref{eq:M_matrix}. There exist a rate $\lambda\in(\max\{\rho(M),\rho_f\},1)$ and constants $\kappa, C_D, C_B, C_\alpha, C_\sigma > 0$ such that for all $t \geq 1$,
\begin{equation}
\begin{aligned}
  \mathbb{E}\big[\|u_t - u_t^{\star}\|\big]
  &\le \kappa\,\lambda^{t-1}
   + C_{D}\sum_{k=1}^{t-1}\lambda^{t-1-k} \|u_{k+1}^{\star} - u_k^{\star}\| \\
  &\quad + C_{B}\sum_{k=1}^{t-1}\lambda^{t-1-k}\,\mathbb{E}\big[\|\theta_{k+1}-\theta_k\|\big]\\
  &\quad + C_{\alpha}\alpha^2 + C_{\sigma}\alpha\sigma_g.
\end{aligned}
\label{eqG:instant_tracking_stoch}
\end{equation}
In particular, if the drifts are uniformly bounded, i.e., 
\begin{equation}
\sup_{t \geq 1} \|u_{t+1}^{\star} - u_t^{\star}\| \leq \bar{\Delta}, \quad \sup_{t \geq 1} \mathbb{E}[\|\theta_{t+1} - \theta_t\|] \leq \bar{B},
\end{equation}
then
\begin{equation}\label{eqG:instant_tracking_stoch_limsup}
\limsup_{t \to \infty} \mathbb{E}[\|u_t - u_t^{\star}\|] \leq \frac{C_D \bar{\Delta} + C_B \bar{B}}{1 - \lambda} + C_\alpha \alpha^2 + C_\sigma\alpha\sigma_g.
\end{equation}
\end{theorem}

Theorem~\ref{thm:tracking_error} shows how S2FO tracks the time-varying optimum in expectation. The asymptotic tracking error depends on the per-step variations $\bar\Delta$ and $\bar B$, the step size $\alpha$, and the stochastic-gradient contribution of order $\alpha\sigma_g$. The term \(C_\sigma\alpha\sigma_g\) quantifies the effect of using a realized stochastic gradient rather than the exact expected gradient.

\subsection{Almost Sure Convergence under Diminishing Step Sizes}\label{sec:almost_sure}
We now consider the asymptotic behavior of S2FO with diminishing step sizes. In this setting, the objective is to characterize convergence to the optimum of a limiting stochastic steady-state problem. This requires the operating condition to settle sufficiently fast relative to the step-size sequence. Specifically, we make the following assumptions.

\begin{assumption}
\label{ass:stepsize}
The step-size sequence \((\alpha_t)_{t\in\mathbb{N}}\) in Algorithm~\ref{alg:ssfo} is nonincreasing and satisfies
\begin{equation}
\alpha_t>0,\qquad
\sum_{t=1}^{\infty}\alpha_t=\infty,
\qquad
\sum_{t=1}^{\infty}\alpha_t^2<\infty.
\end{equation}
\end{assumption}

\begin{assumption}
  \label{ass:slow}
  Let $A_t:=\sum_{s=1}^{t}\alpha_s$ with $(\alpha_s)_{s\in\mathbb{N}}$ being the step sizes of Algorithm~\ref{alg:ssfo}. The stochastic operating condition satisfies
  \begin{equation}
  \sum_{t=1}^{\infty}
  A_t\,\mathbb E\!\left[\|\theta_{t+1}-\theta_t\|\right]
  <\infty.
  \end{equation}
  \end{assumption}

Assumption~\ref{ass:slow} is a weighted summability condition on the drift of the operating condition. It is related to the parameter drift conditions commonly used in non-stationary optimization~\cite{cheung2019learning}. This is stronger than the sublinear-variation condition used in the dynamic regret bound (Remark~\ref{rem:step_size}), which is expected because almost-sure convergence to a fixed limit is a more demanding asymptotic guarantee than sublinear regret. Intuitively, it requires the drift to decay fast enough so that the operating condition approaches a steady regime.
  
  Furthermore, Assumption~\ref{ass:slow} implies that the sequence of distributions \((P_t)\) converges in Wasserstein-1 metric $W_1$. In particular, by the definition of the $W_1$ metric, we have
  \begin{equation}
  W_1(P_{t+1},P_t)
  \le
  \mathbb E\!\left[\|\theta_{t+1}-\theta_t\|\right].
  \end{equation}
  Moreover, since \(A_t\ge \alpha_1>0\), Assumption~\ref{ass:slow} implies that
  \[
  \sum_{t=1}^{\infty}W_1(P_{t+1},P_t)<\infty.
  \]
  Thus, \((P_t)_{t\in\mathbb{N}}\) is a Cauchy sequence in \(W_1\), and we denote its limit by \(P_\infty\in\mathcal P_1(\Theta)\). We then define the limiting cost function and the corresponding minimizer as
  \begin{equation}\label{eq:limit_cost}
  F_\infty(u):=\mathbb E_{\theta\sim P_\infty}[J(u,h(u,\theta))],
  \quad
  u_\infty^\star:=\argmin_{u\in\mathcal U}F_\infty(u).
  \end{equation}
  \noindent\textbf{Example~\ref{ex:wind_decomp} (continued).}
    Continuing with the decomposition \eqref{eq:theta_decomp} introduced in Example~\ref{ex:wind_decomp}, we examine conditions under which Assumption~\ref{ass:slow} holds. By the triangle inequality, the variation satisfies
    \begin{align*}
    \mathbb{E}[\|\theta_{t+1} - \theta_t\|] &\leq \|m_{t+1} - m_t\| + \mathbb{E}[\|\sigma_{t+1}\xi_{t+1} - \sigma_t\xi_t\|]\\
    &\le \|m_{t+1} - m_t\| + (\sigma_{t+1} + \sigma_t)C_\xi.
    \end{align*}
    Therefore, Assumption~\ref{ass:slow} holds if
\begin{equation}
\sum_{t=1}^{\infty}A_t\|m_{t+1}-m_t\|<\infty,
\quad
\sum_{t=1}^{\infty}A_t\sigma_t<\infty.
\end{equation}
For instance, if \(\alpha_t=t^{-a}\) with \(a\in(1/2,1)\), then \(A_t=\mathcal O(t^{1-a})\). Hence the assumption is satisfied when
\[
\|m_{t+1}-m_t\|=\mathcal O(t^{-\beta_m}),
\quad
\sigma_t=\mathcal O(t^{-\beta_\sigma}),
\]
with \(\beta_m>2-a\) and \(\beta_\sigma>2-a\). For another common choice \(\alpha_t=1/t\), we have \(A_t=\mathcal O(\log t)\). The condition then holds if \(\|m_{t+1}-m_t\|=\mathcal O(t^{-\beta_m})\) and 
\(\sigma_t=\mathcal O(t^{-\beta_\sigma})\) with 
\(\beta_m>1\) and \(\beta_\sigma>1\). In both cases, Assumption~\ref{ass:slow} requires the deterministic trend to vary slowly and the stochastic magnitude to decay over time.
  
\begin{lemma}
\label{lem:limiting_cost}
Under Assumption~\ref{ass:slow}, the gradient drift between the stochastic steady-state cost $F_t$ in \eqref{eq:F_t} and its limit $F_\infty$ in \eqref{eq:limit_cost} is summable in the following weighted sense:
\begin{equation}\label{eq:limiting_cost}
\sum_{t=1}^{\infty}
\alpha_t
\sup_{u\in\mathcal U}
\|\nabla F_t(u)-\nabla F_\infty(u)\|
<\infty .
\end{equation}
\end{lemma}

We now show that S2FO converges exactly to the limiting stochastic optimum.
\begin{theorem}
\label{thm:as_convergence}
Suppose Assumptions~\ref{ass:system}--\ref{ass:slow} hold and the limiting cost \(F_\infty\) in \eqref{eq:limit_cost} is $\mu_\infty$-strongly convex on \(\mathcal U\). The iterates of Algorithm~\ref{alg:ssfo} satisfy
\begin{equation}
\lim_{t \to \infty} u_t = u_\infty^{\star} \quad \text{almost surely}.
\end{equation}
\end{theorem}

Theorem~\ref{thm:as_convergence} complements the constant-step-size tracking result in Theorem~\ref{thm:tracking_error}. With a constant step size, the iterates generally remain in a neighborhood of the time-varying optimum, whose size depends on the step size and the stochastic error in the gradient estimate. With diminishing step sizes, the effect of these residual errors shrinks over time, and exact convergence to \(u_\infty^\star\) follows when the operating condition settles sufficiently fast. The tradeoff is that diminishing step sizes typically lead to slower adaptation compared to constant step sizes.

\section{Numerical Simulations}
\label{sec:results}
This section evaluates the proposed S2FO algorithm on a wind farm power maximization problem. We first describe the stochastic wind model used to generate the time-varying inflow conditions. We then present a detailed performance and robustness analysis on a benchmark $2\times 3$ layout under stable wind conditions. We further investigate how the results change as the wind-direction variability increases. Finally, we assess whether the observed trends persist on a subset of the Lillgrund offshore wind farm~\cite{creech2015simulations}.

\subsection{Time-varying wind model}
\label{subsec:wind_model}
In wind energy applications, short-term wind fluctuations are commonly represented by low-order stochastic processes. One widely used model is the Ornstein--Uhlenbeck (OU) process \cite{arenas2020stochastic,zarate2013continuous}. Specifically, we adopt an OU process with a slowly varying mean and fluctuation level, whose continuous-time form is
\begin{equation}
    d\theta_t = \kappa(m_t-\theta_t)\,dt + \sigma_{t}\,dW_t,
    \label{eq:ou}
\end{equation}
where $m_t$ is the slowly varying mean, $\kappa>0$ is the mean-reversion rate, $\sigma_{t}$ is the diffusion coefficient, and $W_t$ is a standard Wiener process. For numerical simulation, we use the corresponding discrete-time AR(1) representation with time-varying parameters:
\begin{equation}
    \theta_{k+1} = m_{k+1} + \rho(\theta_k - m_k) + \epsilon_k, \quad \epsilon_k \sim \mathcal{N}(0,\Sigma_k^2),
    \label{eq:ar1}
\end{equation}
where $\rho\in (0,1)$ is the autoregressive coefficient and $\Sigma_k$ is the noise standard deviation. Equivalently, the inflow in \eqref{eq:ar1} can be decomposed as
\begin{equation}
    \theta_k = m_k + \xi_k,
    \label{eq:decomposition}
\end{equation}
with $m_k$ representing the time-varying deterministic trend and $\xi_k$ representing a zero-mean stochastic fluctuation satisfying $\xi_{k+1} = \rho \xi_k + \epsilon_k$. The AR(1) representation captures the temporal correlation in short-term wind variations, while retaining the mean-plus-fluctuation structure in \eqref{eq:theta_decomp}, as introduced in Example~\ref{ex:wind_decomp}.

\begin{figure}[t]
\centering
\includegraphics[width=0.99\linewidth]{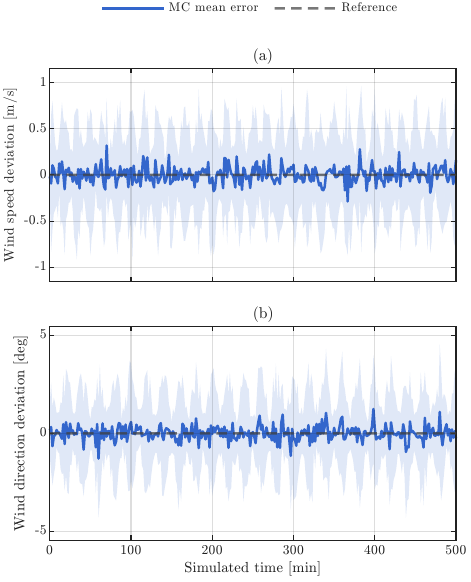}
\caption{Deviation of the Monte Carlo mean of the wind speed and direction from the prescribed reference trend.}
\label{fig:wind_trajectories}
\end{figure}

Fig.~\ref{fig:wind_trajectories} illustrates the typical behavior of this flow model. The figure shows the deviation of the Monte Carlo mean wind speed and direction from their prescribed deterministic trends $m_k$. For both quantities, this deviation remains small and fluctuates around zero without systematic drift, which matches the modeling assumption that $\xi_k$ has zero mean. The model thus produces temporally correlated and non-stationary realizations, which are used in the subsequent controller evaluation.

\subsection{Simulation setup}
\label{sec:setup}

We consider the $2\times 3$ NREL 5\,MW benchmark layout used in \cite{vali2019adjoint}. The turbines are spaced $5D$ in the streamwise direction, with the middle row laterally offset by $0.5D$ relative to the upstream and downstream rows. The wake dynamics are simulated using WFSim~\cite{boersma2018control} on a $2500\times 1500\,\mathrm{m}^2$ domain discretized on a $25\times 15$ staggered grid, with sample time $\Delta t=1$\,s and total horizon $T=3\times 10^4$\,s.

The control inputs are the disk-based thrust coefficient $C_{T,i}'$ and yaw angle $\gamma_i$ of each turbine, constrained to $C_{T,i}'\in[0.4,3.6]$ and $\gamma_i\in[-30^\circ,30^\circ]$. Following \cite{huang2025sequential}, the cost function $J$ in \eqref{eq:F_t} is specified as
\begin{equation}\label{eq:Jsim}
J(C_T',\gamma,y) = \left(\nicefrac{\mathbf{1}_N^\top y}{P^{\mathrm{ref}}} - 1\right)^2 + \frac{\mu}{2}\|C_T'\|^2 + \frac{\mu_\gamma}{2}\|\gamma\|^2,
\end{equation}
with $P^{\mathrm{ref}}=12\,\mathrm{MW}$, $\mu=2.8\times 10^{-4}$, and $\mu_\gamma=2\times 10^{-6}$. The S2FO step sizes are $\alpha=0.3$ for the thrust coefficients and $\alpha_\gamma=4$ for the yaw angles. All controllers are initialized at the greedy steady state, with $C_{T,i}'=2$ and $\gamma_i=0^\circ$.

The wind condition $\theta_t = (v_k,\phi_k)$ collects the inflow speed $v_k$ and direction $\phi_k$, which evolve independently according to the AR(1) model~\eqref{eq:ar1}. The deterministic trends are sinusoidal:
\begin{equation*}
    m_{v,k} = \bar v + A_v\sin(2\pi k/T_v),\quad m_{\phi,k} = \bar\phi + A_\phi\sin(2\pi k/T_\phi).
\end{equation*}
The prescribed fluctuation standard deviations vary as
\begin{equation*}
    \sigma_{v,k} = |\sigma_v + a_v\cos(2\pi k/\tau_v)|,\quad \sigma_{\phi,k} = |\sigma_\phi + a_\phi\cos(2\pi k/\tau_\phi)|.
\end{equation*}
Let $\rho_v$ and $\rho_\phi$ denote the AR(1) coefficients for the wind-speed and wind-direction fluctuations, respectively. The innovation standard deviations in~\eqref{eq:ar1} are chosen as $\Sigma_{v,k}=\sigma_{v,k}\sqrt{1-\rho_v^2}$ and $\Sigma_{\phi,k}=\sigma_{\phi,k}\sqrt{1-\rho_\phi^2}$, so that $\sigma_{u,k}$ and $\sigma_{\phi,k}$ are the stationary standard deviations of the fluctuations $\xi_k$. The nominal parameter values are listed in Table~\ref{tab:wind_params}. They define the low-variability wind condition used in Subsection~\ref{subsec:stable_benchmark}. In Subsections~\ref{subsec:neutral_benchmark} and~\ref{subsec:realistic_layout}, the direction-trend amplitude $A_\phi$ is varied to study increasingly variable wind regimes.

\begin{table}[t]
\centering
\caption{Nominal parameters of the time-varying wind model.}
\label{tab:wind_params}
\begin{tabular}{lcc}
\toprule
Parameter & Speed $v_k$ & Direction $\phi_k$ \\
\midrule
Trend mean & $8\,\mathrm{m/s}$ & $0^\circ$ \\
Trend amplitude & $1.5\,\mathrm{m/s}$ & $5^\circ$ \\
Trend period & $1500\,\mathrm{s}$ & $2000\,\mathrm{s}$ \\
Base std & $0.5\,\mathrm{m/s}$ & $2^\circ$ \\
Std variation & $0.2\,\mathrm{m/s}$ & $1^\circ$ \\
Std period & $1000\,\mathrm{s}$ & $1200\,\mathrm{s}$ \\
AR(1) coefficient & $0.95$ & $0.90$ \\
\bottomrule
\end{tabular}
\end{table}

\subsection{Benchmark layout under low-variability wind conditions}
\label{subsec:stable_benchmark}
 
We begin with the low-variability wind condition, i.e., $A_{\phi} = 5^\circ$. This scenario serves as a controlled reference case for the detailed robustness analysis below. 

\begin{figure}[t]
\centering
\includegraphics[width=0.99\linewidth]{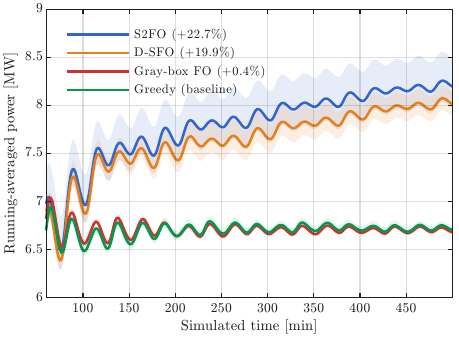}
\caption{Running-averaged farm power on the benchmark $2\times 3$ layout under the low-variability wind condition. Solid lines show the Monte Carlo mean and shaded regions show the $95\%$ confidence interval of the mean.}
 \label{fig:stable_baseline}
\end{figure}

We compare the closed-loop performance of S2FO with three baselines. The first is deterministic SFO (D-SFO), which linearizes the WFSim dynamics around the mean wind condition $m_k$ rather than the realized $\theta_k$. The second is a gray-box FO baseline~\cite{he2026gray}, implemented using data-driven steady-state sensitivity estimation with recursive least squares and a model-free correction. The third is the greedy controller. The comparison is reported in Fig.~\ref{fig:stable_baseline}. S2FO achieves the highest time-averaged farm power, with a gain of $22.7\%$ over the greedy baseline, while D-SFO achieves a gain of $19.9\%$. In contrast, the gray-box FO baseline achieves only a marginal gain of $0.4\%$. A possible explanation is the wake-dominated nature of the wind farm control problem. The effect of a turbine-control change on downstream power is delayed by wake propagation, so a correction based on recent input-output changes may not capture the relevant steady-state wake sensitivity accurately. In addition, the performance gap between S2FO and D-SFO is consistent with Remark~\ref{rem:mean_based_update}: the mean-based sensitivity used in D-SFO does not fully capture the effect of stochastic variations on the expected steady-state performance, while S2FO evaluates the sensitivity at the realized wind condition.

Next, we test how S2FO responds to errors in the wind statistics used for the sensitivity update. We consider two types of error, namely bias in the estimated mean and mismatch in the estimated standard deviation.

\begin{figure}[t]
\centering
\includegraphics[width=0.99\linewidth]{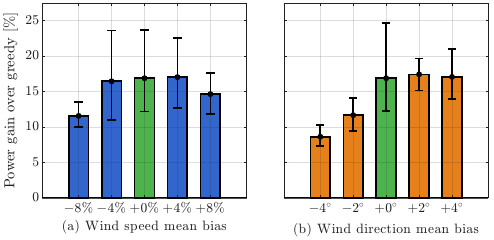}
\caption{Power gain of S2FO over the greedy controller under bias in the assumed mean wind speed and direction. Error bars denote $95\%$ confidence intervals.}
 \label{fig:stable_sensitivity_mean}
\end{figure}

Fig.~\ref{fig:stable_sensitivity_mean} reports the effect of mean-parameter bias in wind speed and wind direction, with error bars showing variability across stochastic inflow realizations. Across all tested biases, S2FO maintains a positive power gain over the greedy baseline. For wind-speed mean bias, the average gain remains close to the unbiased case for moderate biases and decreases mainly at the largest negative bias. Wind-direction mean bias shows a more asymmetric response. Negative biases lead to lower average gains, while positive biases give gains comparable to the unbiased case in the tested range. This difference can be understood from the directional nature of wake interactions. A wind-direction error shifts the expected inflow direction and changes the relative alignment between upstream and downstream turbines, which directly affects the wake pattern used in the sensitivity update. In contrast, a wind-speed bias changes the available power level and wake strength, with a less direct effect on the geometric wake alignment.

\begin{figure}[t]
\centering
\includegraphics[width=0.99\linewidth]{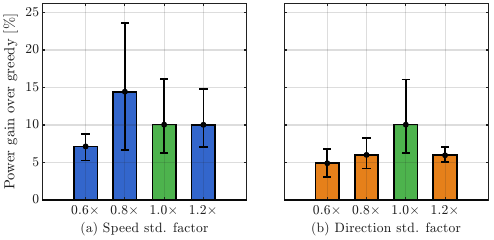}
\caption{Power gain of S2FO over the greedy controller under mismatch in the assumed wind-speed and wind-direction standard deviations. Error bars denote $95\%$ confidence intervals.}
 \label{fig:stable_sensitivity_var}
\end{figure}

Fig.~\ref{fig:stable_sensitivity_var} shows the effect of standard-deviation mismatch. For wind-direction variance, both under-estimation and over-estimation reduce the average gain. This can be explained by the role of direction standard deviation in determining the range of wake-alignment conditions considered in the stochastic sensitivity update. Under-estimating the spread may make the update too specific to a narrow range of wind directions, whereas over-estimating it may average over too broad a range of wake configurations. For wind-speed, the gain also changes with the assumed level, but the tested cases do not show the same clear degradation when deviating from the true value. 

Overall, inaccurate wind statistics generally reduce the achieved power improvement but do not eliminate the benefit of S2FO over the greedy baseline in the tested range. This indicates that the proposed method is robust to moderate statistical mismatch, while also showing that accurate wind-direction statistics are important for obtaining the largest gains. Two mismatched settings, namely the $+2^\circ$ direction bias in Fig.~\ref{fig:stable_sensitivity_mean} and the $0.8\times$ wind-speed standard deviation in Fig.~\ref{fig:stable_sensitivity_var}, yield slightly larger gains than the true-parameter case. A possible explanation is that the biased sensitivity update partially compensates for local linearization or model-approximation errors over the finite simulation horizon.

\subsection{Benchmark layout under increasingly variable wind conditions}
\label{subsec:neutral_benchmark}

Having analyzed the stable reference case in detail, we next examine how the benchmark-layout results change when the wind-direction variability increases. Specifically, we increase the amplitude of the direction trend to $A_\phi = 8^\circ$, corresponding to a moderate-variability wind condition.

\begin{figure}[t]
\centering
\includegraphics[width=0.99\linewidth]{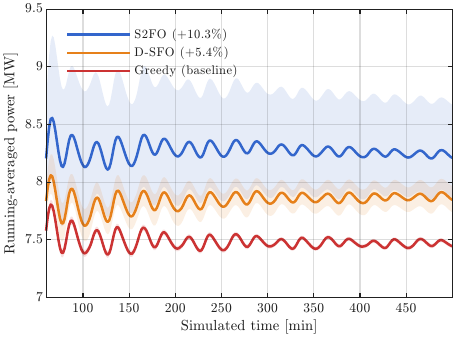}
\caption{Running-averaged farm power on the benchmark $2\times 3$ layout under the moderate-variability wind condition. Solid lines show the Monte Carlo mean and shaded regions show the $95\%$ confidence interval of the mean.}
 \label{fig:neutral_baseline}
\end{figure}

Fig.~\ref{fig:neutral_baseline} shows that, as in the low-variability case, S2FO and D-SFO maintain a higher time-averaged power than the greedy controller throughout the simulation horizon. S2FO achieves a $10.3\%$ gain over the greedy baseline, while D-SFO achieves a $5.4\%$ gain. Both gains are smaller than in the low-variability case, indicating that the improvement margin is reduced under stronger wind-direction variability.

A possible explanation for this reduction lies in the behavior of the greedy baseline itself. As the wind-direction variability increases, the time-averaged greedy power also increases. This increase is related to the wake-alignment pattern. When the inflow direction fluctuates more strongly, downstream turbines spend less time in persistently waked conditions generated by the same upstream turbines. This raises the greedy baseline and leaves a smaller margin for additional improvement by wake-aware feedback optimization. The smaller gain under the moderate-variability condition therefore reflects the stronger greedy baseline at higher wind-direction variability. The same trend is examined on the realistic Lillgrund layout below, where Fig.~\ref{fig:realistic_multi} also shows the increase of the greedy power with amplitude.

\subsection{Realistic wind-farm layout: validation under multiple wind conditions}
\label{subsec:realistic_layout}

We now evaluate the proposed controller on a realistic offshore layout. We consider a \(3\times 4\) subset of the Lillgrund wind farm, whose layout has been widely used as a reference case in wind farm control studies~\cite{van2021sensitivity,creech2015simulations}. The selected subset contains \(N=12\) Siemens SWT-2.3-93 turbines with rotor diameter $D=93\,\mathrm{m}$. The streamwise and lateral spacings between turbines in the subset are $S_x=4.3D$ and $S_y=3.3D$, respectively. The subset is used to preserve realistic wake interactions while keeping the repeated Monte Carlo WFSim simulations computationally tractable.

\begin{figure}[t]
    \centering
    \includegraphics[width=0.99\linewidth]{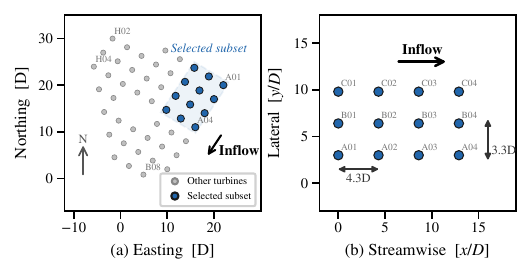}
    \caption{The Lillgrund offshore wind farm and the selected $3\times 4$ subset.}
    \label{fig:lillgrund_layout}
\end{figure}

Fig.~\ref{fig:lillgrund_layout} shows the selected subset in the original Lillgrund layout and its corresponding local coordinate system, in which the nominal inflow direction is aligned with the positive \(x\)-axis.

In WFSim, the simulations use a computational domain of $L_x\times L_y= 2000\times 1000\,\mathrm{m}^2$, discretized with $N_x\times N_y= 36\times 20$ cells. The same stochastic inflow model as in the benchmark study is used, while the wind-direction variability amplitude is varied to study how the achievable gain changes across operating conditions.

\begin{figure}[t]
    \centering
    \includegraphics[width=0.99\linewidth]{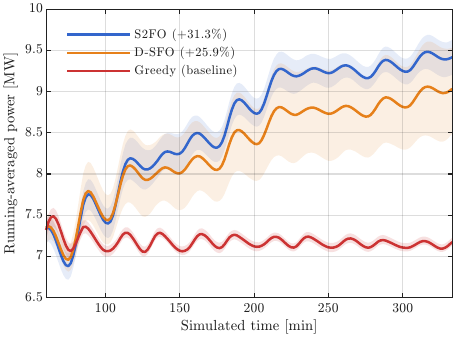}
    \caption{Running-averaged farm power on the Lillgrund $3\times 4$ subset under the low-variability wind condition. Solid lines show the Monte Carlo mean and shaded regions show the $95\%$ confidence interval of the mean.}
    \label{fig:realistic_stable_threeway}
\end{figure}

Fig.~\ref{fig:realistic_stable_threeway} gives the baseline comparison under the low-variability wind condition. The time-averaged power trajectories show that both SFO versions outperform the greedy controller after the initial transient. D-SFO improves the greedy baseline by \(25.9\%\), while S2FO achieves a larger gain of \(31.3\%\). Thus, the stochastic advantage observed in the benchmark layout persists in the Lillgrund subset.

\begin{figure}[t]
    \centering
    \includegraphics[width=0.99\linewidth]{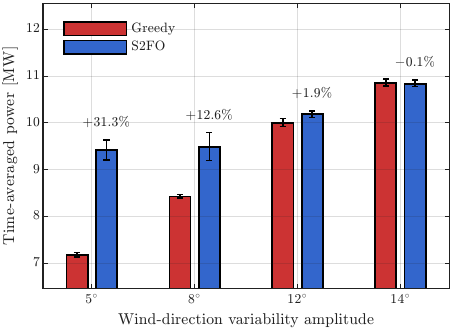}
    \caption{Time-averaged power of the greedy controller and S2FO on the Lillgrund $3\times 4$ subset under different wind-direction variability amplitudes. Error bars denote $95\%$ confidence intervals, and the percentages above the bars indicate the gain of S2FO over the greedy baseline.}
    \label{fig:realistic_multi}
\end{figure}

Fig.~\ref{fig:realistic_multi} summarizes the comparison between the greedy controller and S2FO under different wind-direction variability amplitudes. The improvement over greedy control decreases as the amplitude increases. The gain is \(31.3\%\) at \(5^\circ\), \(12.6\%\) at \(8^\circ\), \(1.9\%\) at \(12^\circ\), and approximately \(-0.1\%\) at \(14^\circ\). Hence, S2FO provides a clear benefit under stable and moderately varying wind conditions, but the gain becomes negligible under the strongest variability tested here.

This trend is consistent with the benchmark-layout results in Section~\ref{subsec:neutral_benchmark}. In both layouts, stronger wind-direction variability weakens persistent wake alignment, raises the greedy baseline, and leaves a smaller margin for wake-aware control. Therefore, the benefit of S2FO is more evident when persistent wake interactions leave sufficient room for coordinated control.

\section{Conclusion}
This paper proposes stochastic sequential feedback optimization for steady-state optimization of nonlinear dynamical systems in random operating environments. The method evaluates the sensitivity update at realized operating conditions rather than only at a nominal condition. The analysis shows that the performance of S2FO is governed by the movement of the stochastic optimizer, the plant transient relative to the steady-state map, and the variation of the operating condition. The finite-horizon bounds quantify how these terms affect regret and tracking, while the result with diminishing step sizes shows that exact convergence can be recovered when the operating condition stabilizes sufficiently fast.

Wind farm simulations indicate that incorporating stochastic operating-condition information in the sensitivity update can improve average farm power compared with greedy control and mean-based deterministic SFO. The benefit is most visible when persistent wake interactions leave room for coordinated control. When wind-direction variability becomes stronger, this margin decreases, which in turn reduces the improvement over the greedy baseline.

Future work will consider online estimation of the operating-condition distribution, less conservative sensitivity approximations, and extensions to larger wind farms with more detailed atmospheric and actuator models.

\section*{Appendix}
\subsection{Proofs in Section~\ref{sec:pre_lemmas}}
\subsubsection*{Proof of Lemma~\ref{lem:sensitivity_lip}}
For notational convenience, denote $A_i := \nabla_x f(x_i, u_i, \theta_i)$ and $B_i := \nabla_u f(x_i, u_i, \theta_i)$ for $i = 1, 2$. Then
\begin{align*}
&\|\widetilde{\nabla h}(x_1, u_1, \theta_1) - \widetilde{\nabla h}(x_2, u_2, \theta_2)\| \notag \\
&= \|(I - A_1)^{-1} B_1 - (I - A_2)^{-1} B_2\| \notag \\
&= \|(I - A_1)^{-1}(B_1 - B_2) + ((I - A_1)^{-1} - (I - A_2)^{-1}) B_2\| \notag \\
&\leq \|(I - A_1)^{-1}\| \|B_1 - B_2\|\notag\\
&\quad + \|(I - A_1)^{-1} - (I - A_2)^{-1}\| \|B_2\|.
\end{align*}
By Assumption~\ref{ass:system}(i), $\|A_1\| \leq \rho_f < 1$, so
\begin{equation*}
\|(I - A_1)^{-1}\| \leq \sum_{k=0}^{\infty} \|A_1\|^k \leq \frac{1}{1 - \rho_f}.
\end{equation*}
For the difference of inverses, using $A^{-1} - B^{-1} = A^{-1}(B - A)B^{-1}$ yields
\begin{align*}
&\quad\|(I - A_1)^{-1} - (I - A_2)^{-1}\|\notag\\
 &= \|(I - A_1)^{-1}(A_1 - A_2)(I - A_2)^{-1}\| \notag \\
&\leq \frac{1}{(1-\rho_f)^2} \|A_1 - A_2\| \notag \\
&\leq \frac{L_{f,x}}{(1-\rho_f)^2} (\|x_1 - x_2\| + \|u_1 - u_2\| + \|\theta_1 - \theta_2\|),
\end{align*}
where we used the Lipschitz continuity of $\nabla_x f$ from Assumption~\ref{ass:system}(ii). Similarly, $\|B_1 - B_2\| \leq L_{f,u}(\|x_1 - x_2\| + \|u_1 - u_2\| + \|\theta_1 - \theta_2\|)$ and $\|B_2\| \leq G_u^f$. Combining these bounds gives the result with $C = \frac{L_{f,u}}{1-\rho_f} + \frac{G_u^f L_{f,x}}{(1-\rho_f)^2} = \frac{(1-\rho_f)L_{f,u} + G_u^f L_{f,x}}{(1-\rho_f)^2}$. \hfill $\blacksquare$

\subsubsection*{Proof of Lemma~\ref{lem:output_error}}
By the triangle inequality,
\begin{align}
e_{t+1} &= \|y_{t+1} - h(u_{t+1}, \theta_{t+1})\| \notag \\
&\leq \|y_{t+1} - h(u_t, \theta_t)\| + \|h(u_t, \theta_t) - h(u_{t+1}, \theta_{t+1})\|.
\label{eq:e_split}
\end{align}
For the first term, using the dynamics $y_{t+1} = f(y_t, u_t, \theta_t)$ and the fixed-point relation $h(u_t, \theta_t) = f(h(u_t, \theta_t), u_t, \theta_t)$, we have
\begin{align}
\|y_{t+1} - h(u_t, \theta_t)\| &= \|f(y_t, u_t, \theta_t) - f(h(u_t, \theta_t), u_t, \theta_t)\| \notag \\
&\leq \rho_f \|y_t - h(u_t, \theta_t)\| = \rho_f e_t,
\label{eq:e_term1}
\end{align}
where we used Assumption~\ref{ass:system}(i). For the second term, by Lipschitz continuity of $h$ in $(u, \theta)$ from Assumption~\ref{ass:system}(iii),
\begin{equation}
\|h(u_t, \theta_t) - h(u_{t+1}, \theta_{t+1})\| \leq L_{h,u} \|u_{t+1} - u_t\| + L_{h,\theta} \|\theta_{t+1} - \theta_t\|.
\label{eq:e_term2}
\end{equation}
Setting $C_u := L_{h,u}$ and $C_\theta := L_{h,\theta}$ and substituting \eqref{eq:e_term1}--\eqref{eq:e_term2} into \eqref{eq:e_split} yields \eqref{eq:e_recursion}. Unrolling the recursion gives
\begin{equation*}
e_t \leq \rho_f^{t-1} e_1 + \sum_{k=1}^{t-1} \rho_f^{t-1-k} \big(C_u \|u_{k+1} - u_k\| + C_\theta \|\theta_{k+1} - \theta_k\|\big).
\end{equation*}
Summing over $t = 1, \ldots, T$, we get
\begin{align*}
\sum_{t=1}^{T} e_t &\leq \sum_{t=1}^{T} \rho_f^{t-1} e_1 + \sum_{t=1}^{T}\sum_{k=1}^{t-1} \rho_f^{t-1-k} C_u \|u_{k+1} - u_k\|\notag\\
&\quad + \sum_{t=1}^{T}\sum_{k=1}^{t-1}\rho_f^{t-1-k}C_\theta \|\theta_{k+1} - \theta_k\| \notag \\
&\leq \frac{e_1}{1-\rho_f} + \sum_{k=1}^{T-1} C_u \|u_{k+1} - u_k\| \sum_{j=0}^{\infty} \rho_f^j\notag\\
&\quad + \sum_{k=1}^{T-1}C_\theta \|\theta_{k+1} - \theta_k\| \sum_{j=0}^{\infty} \rho_f^j \notag \\
&= \frac{e_1}{1-\rho_f} + \frac{C_u}{1-\rho_f} \sum_{k=1}^{T-1} \|u_{k+1} - u_k\|\notag\\
&\quad + \frac{C_\theta}{1-\rho_f} \sum_{k=1}^{T-1} \|\theta_{k+1} - \theta_k\|,
\end{align*}
where the second inequality follows from exchanging the order of summation and bounding the geometric series. \hfill $\blacksquare$

\subsubsection*{Proof of Lemma~\ref{lem:grad_error}}
We first prove \eqref{eq:grad_mismatch_bound}. By definition,
\begin{align}
&\quad \nabla \tilde{J}(u_t, \theta_t) - g_t\notag\\
&= \big(\nabla_u J(u_t,h(u_t,\theta_t)) - \nabla_u J(u_t,y_t)\big)\label{lemF:grad_decomp}\\
&\quad + \nabla h(u_t,\theta_t)^{\top}\big(\nabla_y J(u_t,h(u_t,\theta_t)) - \nabla_y J(u_t,y_t)\big)\notag\\
&\quad + \big(\nabla h(u_t,\theta_t) - \widetilde{\nabla h}(y_t,u_t,\theta_t)\big)^{\top}\nabla_y J(u_t,y_t)\notag\\
&=: A_t + B_t + C_t.\notag
\end{align}
By Assumption~\ref{ass:cost}(i), $\nabla_u J$ is Lipschitz in $y$ with constant $L_{J,u}$, hence
\begin{equation}
\label{lemF:A_bound}
\|A_t\|
= \|\nabla_u J(u_t,h(u_t,\theta_t)) - \nabla_u J(u_t,y_t)\|
\le L_{J,u}e_t.
\end{equation}
Furthermore, we have
\begin{align*}
\|B_t\|
&= \big\|\nabla h(u_t,\theta_t)^{\top}\big(\nabla_y J(u_t,h(u_t,\theta_t)) - \nabla_y J(u_t,y_t)\big)\big\|\notag\\
&\le \|\nabla h(u_t,\theta_t)\|\cdot\|\nabla_y J(u_t,h(u_t,\theta_t)) - \nabla_y J(u_t,y_t)\|.
\label{lemF:B_step1}
\end{align*}
By Assumption~\ref{ass:system}, there exists $\hat{L}_h>0$ such that
\begin{equation*}
\label{lemF:nabla_h_bound}
\|\nabla h(u_t,\theta_t)\|\le \hat{L}_h.
\end{equation*}
Moreover, the Lipschitz continuity of $\nabla_y J$ in $y$ gives
\[
\|\nabla_y J(u_t,h(u_t,\theta_t)) - \nabla_y J(u_t,y_t)\|
\le L_{J,y}e_t.
\]
Combining the two bounds yields
\begin{equation}
\label{lemF:B_bound}
\|B_t\|\le \hat{L}_hL_{J,y}e_t.
\end{equation}

By differentiating the fixed-point relation
$h(u,\theta)=f(h(u,\theta),u,\theta)$ with respect to $u$,
\begin{align*}
\nabla h(u,\theta)
&=
\bigl(I-\nabla_x f(h(u,\theta),u,\theta)\bigr)^{-1}
\nabla_u f(h(u,\theta),u,\theta)\\
&=
\widetilde{\nabla h}(h(u,\theta),u,\theta),
\end{align*}
where the inverse exists by Assumption~\ref{ass:system}(i). Therefore, Lemma~\ref{lem:sensitivity_lip} gives
\begin{align*}
&\quad\|\nabla h(u_t,\theta_t) - \widetilde{\nabla h}(y_t,u_t,\theta_t)\|\notag\\
&= \|\widetilde{\nabla h}(h(u_t,\theta_t),u_t,\theta_t) - \widetilde{\nabla h}(y_t,u_t,\theta_t)\|\notag\\
&\le C\,\|h(u_t,\theta_t) - y_t\| = Ce_t.
\label{lemF:C_step1}
\end{align*}
Together with \(\|\nabla_yJ(u_t,y_t)\|\le G_y^J\), this implies 
\begin{equation}
\label{lemF:C_bound}
\|C_t\|
\le G_y^J\,\|\nabla h(u_t,\theta_t) - \widetilde{\nabla h}(y_t,u_t,\theta_t)\|
\le G_y^J Ce_t.
\end{equation}
Substituting \eqref{lemF:A_bound}, \eqref{lemF:B_bound} and \eqref{lemF:C_bound} into \eqref{lemF:grad_decomp}, we obtain
\begin{equation*}\label{eq:term1}
\|\nabla \tilde{J}(u_t, \theta_t) - g_t\|
\le (L_{J,u} + \hat{L}_hL_{J,y} + G_y^J C)e_t.
\end{equation*}

It remains to prove \eqref{eq:cond_unbiased}. First, by the chain rule and the uniform bounds in Assumptions~\ref{ass:system}--\ref{ass:cost}, $\|\nabla \tilde J(u,\theta)\|
\leq
G_u^J+\hat L_h G_y^J
$ holds uniformly in $(u,\theta)$. Hence, differentiation and expectation can be interchanged by the dominated convergence theorem, so that
\[
\nabla F_t(u)=\mathbb{E}_{\theta_t\sim P_t}
\left[\nabla\tilde J(u,\theta_t)\right].
\]
Since $u_t$ is $\mathcal F_{t-1}$-measurable and $\theta_t$ is drawn from $P_t$ independently of $\mathcal F_{t-1}$, it follows that
\begin{align*}
\mathbb{E}\left[
\nabla \tilde J(u_t,\theta_t)
\mid \mathcal F_{t-1}
\right]
&=
\mathbb{E}_{\theta_t\sim P_t}
\left[
\nabla \tilde J(u_t,\theta_t)
\right]
=
\nabla F_t(u_t).
\end{align*}
This proves the lemma. \hfill $\blacksquare$

\subsection{Proofs in Section~\ref{sec:regret}}
\subsubsection*{Proof of Theorem~\ref{thm:dynamic_regret}}
By convexity of $F_t$,
\begin{equation*}
F_t(u_t) - F_t(u_t^{\star}) \leq \nabla F_t(u_t)^\top(u_t - u_t^{\star}).
\end{equation*}
Using Lemma~\ref{lem:grad_error}, we decompose
\begin{align}
\nabla F_t(u_t)^\top(u_t-u_t^\star)
&=
g_t^\top(u_t-u_t^\star) \notag\\
&\quad
+\big(\nabla \tilde J(u_t,\theta_t)-g_t\big)^\top(u_t-u_t^\star) \notag\\
&\quad
+\big(\nabla F_t(u_t)-\nabla \tilde J(u_t,\theta_t)\big)^\top(u_t-u_t^\star).\label{eq:inst_regret_decomp}
\end{align}
The second term bounded by $2R\bar L e_t$. Define $\zeta_t := \nabla \tilde{J}(u_t, \theta_t) - \nabla F_t(u_t)$. Then, by \eqref{eq:cond_unbiased}, \(\mathbb E[\zeta_t\mid\mathcal F_{t-1}]=0\). Since \(u_t\) and \(u_t^\star\) are \(\mathcal F_{t-1}\)-measurable,
\begin{equation}\label{eq:independent}
\mathbb E[\zeta_t^\top(u_t-u_t^\star)\mid\mathcal F_{t-1}]=0.
\end{equation}

By the nonexpansiveness of the projection,
\begin{align*}
\|u_{t+1} - u_t^{\star}\|^2 &\leq \|u_t - \alpha g_t - u_t^{\star}\|^2 \notag \\
&= \|u_t - u_t^{\star}\|^2 - 2\alpha g_t^\top(u_t - u_t^{\star}) + \alpha^2 \|g_t\|^2.
\end{align*}
Rearranging the terms gives
\begin{align*}
  g_t^\top(u_t-u_t^\star)
  &\le
  \frac{1}{2\alpha}
  \left(
  \|u_t-u_t^\star\|^2-\|u_{t+1}-u_t^\star\|^2
  \right)
  +\frac{\alpha}{2}\|g_t\|^2\\
  &= \frac{1}{2\alpha}
  \left(
  \|u_t-u_t^\star\|^2-\|u_{t+1}-u_{t+1}^\star\|^2
  \right)
  +\frac{\alpha}{2}\|g_t\|^2\\
  &\quad + \frac{1}{2\alpha}\left(\|u_{t+1}-u_{t+1}^\star\|^2 - \|u_{t+1}-u_t^\star\|^2\right)
  \end{align*}
  For the last term, we have
  \begin{align*}
  &\quad\|u_{t+1}-u_{t+1}^\star\|^2 - \|u_{t+1}-u_t^\star\|^2\\
  &= (2u_{t+1} - u_{t+1}^\star - u_t^\star)^\top(u_t^\star - u_{t+1}^\star)\\
  &\le 4R\|u_{t+1}^\star - u_t^\star\|.
  \end{align*}
  Therefore, we get
\begin{align}\label{eq:regret_term1}
g_t^\top(u_t - u_t^{\star}) &\leq \frac{1}{2\alpha}\big(\|u_t - u_t^{\star}\|^2 - \|u_{t+1} - u_{t+1}^{\star}\|^2\big)\notag\\
&\quad + \frac{4R}{2\alpha}\|u_{t+1}^{\star} - u_t^{\star}\| + \frac{\alpha}{2}\|g_t\|^2.
\end{align}
Summing \eqref{eq:inst_regret_decomp} over $t = 1, \ldots, T$ yields
\begin{align*}
\mathrm{reg}(T) &\leq \sum_{t=1}^{T} g_t^\top(u_t - u_t^{\star}) + 2R \bar{L} \sum_{t=1}^{T} e_t - \sum_{t=1}^{T} \zeta_t^\top(u_t - u_t^{\star}).
\end{align*}
Taking expectations and using \eqref{eq:independent},
\begin{align*}
\mathbb{E}[\mathrm{reg}(T)] &\leq \mathbb{E}\left[\sum_{t=1}^{T} g_t^\top(u_t - u_t^{\star})\right] + 2R \bar{L} \mathbb{E}\left[\sum_{t=1}^{T} e_t\right].
\end{align*}
Under Assumptions~\ref{ass:system}--\ref{ass:cost} and the boundedness of $\mathcal{U}$, the gradient estimate $g_t$ is uniformly bounded. We denote by $G$ a uniform upper bound on $\|g_t\|$. Then, we can use \eqref{eq:regret_term1} to bound the first term and use Lemma~\ref{lem:output_error} to bound the second term. Combining these bounds yields \eqref{eq:regret_bound}. \hfill $\blacksquare$

\subsection{Proofs in Section~\ref{sec:instan_tracking}}
\subsubsection*{Proof of Theorem~\ref{thm:tracking_error}}
Define $y_t^{\star}:= h(u_t^{\star}, \theta_t)$ and the instantaneous tracking errors
\[
  e_t^u := \|u_t - u_t^{\star}\|,\qquad
  e_t^y := \|y_t - y_t^{\star}\|.
\]
From the fixed-point relation, we have
\[
u_t^\star=\Pi_{\mathcal U}[u_t^\star-\alpha\nabla F_t(u_t^\star)].
\]
By non-expansiveness of the projection,
\[
\begin{aligned}
\|u_{t+1}-u_t^\star\|
&\le
\left\|
u_t-u_t^\star-\alpha\big(g_t-\nabla F_t(u_t^\star)\big)
\right\| \\
&\le
\left\|
u_t-u_t^\star-\alpha\big(g_t-\nabla\tilde J(u_t^\star,\theta_t)\big)
\right\|
+\alpha\|\zeta_t^\star\|,
\end{aligned}
\]
where $\zeta_t^\star:= \nabla\tilde J(u_t^\star,\theta_t)-\nabla F_t(u_t^\star)$ denotes the stochastic realization error. The first term has the same form as in the deterministic sequential-linearization analysis. Hence, by the contraction estimate in~\cite[Theorem~3.6]{huang2025sequential}, with constants \(C_1,C_2\) defined above and an additional constant $C_3 > 0$,
\begin{align*}
&\quad\left\|
u_t-u_t^\star-\alpha\big(g_t-\nabla\tilde J(u_t^\star,\theta_t)\big)
\right\|\\
&\le
\left(\sqrt{1-2\alpha\mu_J+\alpha^2L_{J,u}^2}+\alpha C_1\right)e_t^u
+\alpha C_2e_t^y + \alpha C_3e_t,
\end{align*}
where $e_t := \|y_t - h(u_t,\theta_t)\|$ is the output mismatch. Therefore,
\begin{align*}
  e_{t+1}^u &= \|u_{t+1} - u_t^\star + u_t^\star - u_{t+1}^\star\|\notag\\
  &\le m_{11} e_t^u + m_{12} e_t^y + \alpha C_3 e_t + \alpha\|\zeta_t^\star\| + \|u_{t+1}^{\star}-u_t^{\star}\|,
\end{align*}
with $m_{11} := \sqrt{1 - 2\alpha\mu_J + \alpha^2L_{J,u}^2} + \alpha C_1$ and $m_{12} := \alpha C_2$. By Assumption~\ref{ass:noise} and Jensen's inequality, $E[\|\zeta_t^\star\|]\leq \sigma_g$. Hence, defining
\[\bar e_t^u:=\mathbb E[e_t^u],
\qquad
\bar e_t^y:=\mathbb E[e_t^y],
\qquad
\bar e_t:=\mathbb E[e_t],\]
we obtain
\begin{align}
  \bar e_{t+1}^u
  &\le m_{11} \bar e_t^u + m_{12} \bar e_t^y + \alpha C_3 \bar e_t + \alpha\sigma_g + \|u_{t+1}^{\star}-u_t^{\star}\|.
  \label{eqG:u_err_rec}
\end{align}

The evolution of the output error follows from the system dynamics and the steady-state relation for $y_t^{\star}$.
Specifically,
\begin{align*}
  e_{t+1}^y
  &= \|f(y_t,u_t,\theta_t) - f(y_{t+1}^{\star},u_{t+1}^{\star},\theta_{t+1})\|\\
  &\le \|f(y_t,u_t,\theta_t) - f(y_t^{\star},u_t^{\star},\theta_t)\|\\
  &\quad + \|f(y_t^{\star},u_t^{\star},\theta_t) - f(y_{t+1}^{\star},u_{t+1}^{\star},\theta_{t+1})\|.
\end{align*}
Using Assumption~\ref{ass:system}(i) and the Lipschitz continuity of $f$ in $(u,\theta)$, we derive
\[
  e_{t+1}^y \le m_{21} e_t^u + m_{22} e_t^y + C_4\|\theta_{t+1}-\theta_t\| + C_5\|u_{t+1}^{\star} - u_t^\star\|,
\]
with $m_{21}:=G_u^f$, $m_{22}:=\rho_f$, $C_4:=L_{f,\theta} + \rho_f L_{h,\theta}$, and $C_5:= G_u^f+ \rho_f L_{h,u}$. Taking expectations yields
\begin{equation}\label{eqG:y_err_rec}
\bar e_{t+1}^y
\leq
m_{21}\bar e_t^u+m_{22}\bar e_t^y
+C_4B_t+C_5D_t,\end{equation}
where $D_t:=\|u_{t+1}^\star-u_t^\star\|$ and $B_t:=\mathbb E[\|\theta_{t+1}-\theta_t\|]$. 

Let $\bar E_t := \begin{bmatrix} \bar e_t^u \\ \bar e_t^y \end{bmatrix}$. Stacking \eqref{eqG:u_err_rec} and \eqref{eqG:y_err_rec}, we obtain
\begin{equation}
  \bar E_{t+1} \le M \bar E_t + b_t,
  \label{eqG:E_rec_mov}
\end{equation}
where $M$ is defined in \eqref{eq:M_matrix} and $b_t:=\begin{bmatrix}
            \alpha C_3\bar e_t + \alpha\sigma_g + D_t\\[2pt]
            C_4B_t + C_5D_t
         \end{bmatrix}$. The assumption $\rho(M)<1$ guarantees that the homogeneous system $\bar E_{t+1} = M \bar E_t$ is exponentially stable. By Lemma~\ref{lem:output_error} and $\|u_{j+1}-u_j\|\le \alpha G$, 
\begin{equation}
  \bar e_t \le \rho_f^{t-1}\bar e_1 + \frac{C_u G}{1-\rho_f} \alpha + C_{\theta}\sum_{j=1}^{t-1}\rho_f^{t-1-j}B_j.
  \label{eqG:e_t_bound_rec}
\end{equation}
Substituting \eqref{eqG:e_t_bound_rec} into $b_t$, we obtain
\begin{align}
  \|b_t\|
  &\le (1+C_5)D_t
     + C_4B_t + \alpha\sigma_g + \alpha C_3C_{\theta}\sum_{j=1}^{t-1}\rho_f^{t-1-j}B_j\notag\\
     &\quad+ \frac{C_3C_uG}{1-\rho_f}\alpha^2 + \alpha C_3\rho_f^{t-1}\bar e_1.
  \label{eqG:b_t_bound}
\end{align}
Iterating \eqref{eqG:E_rec_mov} gives
\[
  \bar E_t \le M^{t-1} \bar E_1 + \sum_{k=1}^{t-1} M^{t-1-k} b_k.
\]
Since \(\rho(M)<1\), for any \(\lambda\in(\rho(M),1)\) there exists \(c_M>0\) such that $\|M^{t-1}\|\le c_M\,\lambda^{t-1}$. Specifically, we choose \(\lambda\in(\max\{\rho(M),\rho_f\},1)\). Then
\[
  \|\bar E_t\|
  \le c_M\,\lambda^{t-1}\|\bar E_1\|
   + c_M\sum_{k=1}^{t-1}\lambda^{t-1-k}\|b_k\|.
\]
Combining with \eqref{eqG:b_t_bound}, we get
\begin{equation}
\begin{aligned}
  \|\bar E_t\|
  &\le c_M\|\bar E_1\|\,\lambda^{t-1}
   + c_M(1+C_5)\sum_{k=1}^{t-1}\lambda^{t-1-k} D_k\\
   &\quad + c_M\alpha C_3C_\theta\sum_{k=1}^{t-1}\lambda^{t-1-k}\sum_{j=1}^{k-1}\rho_f^{k-1-j}B_j\\
  &\quad + c_MC_4\sum_{k=1}^{t-1}\lambda^{t-1-k}B_k\\
  &\quad + \frac{c_M}{1-\lambda}\left(\frac{C_3C_uG}{1-\rho_f}\alpha^2 + \alpha\sigma_g + \alpha C_3\bar e_1\right).
\end{aligned}
\label{eqG:E_t_path}
\end{equation}
By exchanging the order of summation, 
\begin{align}\sum_{k=1}^{t-1}
  \lambda^{t-1-k}
  \sum_{j=1}^{k-1}\rho_f^{k-1-j}
  B_j &=
  \sum_{j=1}^{t-2}
  B_j
  \sum_{k=j+1}^{t-1}
  \lambda^{t-1-k}\rho_f^{k-1-j}\notag\\
  &= \sum_{j=1}^{t-2}
  B_j\sum_{r=0}^{t-2-j}\lambda^{t-2-j-r}\rho_f^r\notag\\
  &\le \frac{1}{\lambda - \rho_f}\sum_{j=1}^{t-2}\lambda^{t-1-j}B_j.\label{eq:double_convolution}
  \end{align}
Since $\bar e_t^u \le \|\bar E_t\|$, inequalities \eqref{eqG:E_t_path} and \eqref{eq:double_convolution} imply \eqref{eqG:instant_tracking_stoch} with $\kappa:=c_M\|\bar E_1\|$, $C_D := c_M(1+C_5)$, $C_B := c_MC_4 + \frac{\alpha C_3C_\theta}{\lambda-\rho_f}$, $C_\alpha := \frac{c_M}{1-\lambda}\left(\frac{C_3C_uG}{1-\rho_f}\right)$, and $C_\sigma := \frac{c_M}{1-\lambda}$.  
If $\sup_t D_t\le \bar \Delta$ and $\sup_t B_t\le \bar B$, we have
\[
  \sum_{k=1}^{t-1}\lambda^{t-1-k} D_k
  \le \bar \Delta \sum_{j=0}^{\infty}\lambda^j
  = \frac{\bar \Delta}{1-\lambda},
\]
and similarly for the parameter variation term. Taking $\limsup$ in \eqref{eqG:instant_tracking_stoch}, we obtain \eqref{eqG:instant_tracking_stoch_limsup}. \hfill $\blacksquare$

\subsection{Proofs in Section~\ref{sec:almost_sure}}
\subsubsection*{Proof of Lemma~\ref{lem:limiting_cost}}
By the Lipschitz continuity of \(\nabla\tilde J(u,\theta)\) in \(\theta\) and Kantorovich-Rubinstein duality of Wasserstein-1 metric \cite[Lemma~2.5]{bianchin2023online}, there exists \(L_{\nabla F,\theta}>0\) such that
\begin{equation}\label{eq:KR-duality}
\|\nabla F_t(u) - \nabla F_\infty(u)\| \leq L_{\nabla F,\theta} W_1(P_t, P_\infty), \qquad
\forall u\in\mathcal U.
\end{equation}
Moreover, by the triangle inequality and the definition of \(P_\infty\),
\begin{equation*}
W_1(P_t,P_\infty)
\le
\sum_{k=t}^{\infty}W_1(P_{k+1},P_k)
\le
\sum_{k=t}^{\infty}
\mathbb E[\|\theta_{k+1}-\theta_k\|].
\end{equation*}
Therefore,
\begin{align*}
\sum_{t=1}^{\infty}\alpha_t W_1(P_t,P_\infty)
&\le
\sum_{t=1}^{\infty}\alpha_t
\sum_{k=t}^{\infty}
\mathbb E[\|\theta_{k+1}-\theta_k\|] \notag\\
&=
\sum_{k=1}^{\infty}
\left(\sum_{t=1}^{k}\alpha_t\right)
\mathbb E[\|\theta_{k+1}-\theta_k\|] \notag\\
&=
\sum_{k=1}^{\infty}
A_k\,\mathbb E[\|\theta_{k+1}-\theta_k\|]
<\infty,
\end{align*}
where the last inequality follows from Assumption~\ref{ass:slow}. Combining this with \eqref{eq:KR-duality} yields \eqref{eq:limiting_cost}. \hfill $\blacksquare$
\subsubsection*{Proof of Theorem~\ref{thm:as_convergence}}
From Lemma~\ref{lem:output_error},
\begin{equation*}
e_{t+1} \leq \rho_f e_t + C_u \|u_{t+1} - u_t\| + C_\theta \|\theta_{t+1} - \theta_t\|.
\end{equation*}
Using the gradient bound $\|g_t\| \leq G$ and the update rule, we have $\|u_{t+1} - u_t\| \leq \alpha_t G$. Hence,
\begin{equation*}
e_{t+1} \leq \rho_f e_t + C_u G \alpha_t + C_\theta \|\theta_{t+1} - \theta_t\|.
\end{equation*}
Unrolling this recursion gives
\begin{equation}
  e_t
  \le
  \rho_f^{t-1}e_1
  +
  \sum_{k=1}^{t-1}\rho_f^{t-1-k}
  \left(
  C_uG\alpha_k
  +
  C_\theta\|\theta_{k+1}-\theta_k\|
  \right).
  \label{eq:as_e_unroll}
  \end{equation}
  Taking expectations in \eqref{eq:as_e_unroll}, multiplying by \(\alpha_t\), and summing over \(t\), we obtain
  \begin{align*}
  \sum_{t=1}^{\infty}\alpha_t\mathbb E[e_t]
  &\le
  e_1\sum_{t=1}^{\infty}\alpha_t\rho_f^{t-1} + \sum_{t=1}^{\infty}
  \alpha_t
  \sum_{k=1}^{t-1}\rho_f^{t-1-k} C_uG\alpha_k\notag\\
  &\quad
  +
  \sum_{t=1}^{\infty}
  \alpha_t
  \sum_{k=1}^{t-1}\rho_f^{t-1-k}
  C_\theta\mathbb E[\|\theta_{k+1}-\theta_k\|] \notag\\
  &=
  e_1\sum_{t=1}^{\infty}\alpha_t\rho_f^{t-1}+\sum_{k=1}^{\infty}C_uG\alpha_k\sum_{t=k+1}^{\infty}
  \alpha_t\rho_f^{t-1-k} \notag\\
  &\quad +
  C_\theta\sum_{k=1}^{\infty}\mathbb E[\|\theta_{k+1}-\theta_k\|]
  \sum_{t=k+1}^{\infty}
  \alpha_t\rho_f^{t-1-k}.
  \end{align*}
Since \((\alpha_t)\) is nonincreasing, \(\alpha_t\le \alpha_k\) for \(t\ge k+1\). Therefore,
\begin{equation*}
\sum_{t=k+1}^{\infty}
\alpha_t\rho_f^{t-1-k}
\le
\frac{\alpha_k}{1-\rho_f}.
\end{equation*}
It follows that
\begin{align}
\sum_{t=1}^{\infty}\alpha_t \mathbb E[e_t]
&\le
\frac{\alpha_1e_1}{1-\rho_f}
+
\frac{C_uG}{1-\rho_f}\sum_{k=1}^{\infty}\alpha_k^2\notag\\
&\quad
+
\frac{C_\theta}{1-\rho_f}
\sum_{k=1}^{\infty}A_k\mathbb E[\|\theta_{k+1}-\theta_k\|] < \infty, \label{eq:weighted_e_L1}
\end{align}
where we used $\alpha_k\le A_k$, together with Assumptions~\ref{ass:stepsize} and~\ref{ass:slow}.

We now rewrite the update as a stochastic approximation to the limiting problem. Define $\zeta_t
:=
\nabla\tilde J(u_t,\theta_t)-\nabla F_t(u_t)$ and $r_t
:=
g_t-\nabla\tilde J(u_t,\theta_t)
+
\nabla F_t(u_t)-\nabla F_\infty(u_t)$.
Then
\[
g_t=\nabla F_\infty(u_t)+r_t+\zeta_t,
\]
and the update can be written as
\[
u_{t+1}
=
\Pi_{\mathcal U}
\left[
u_t-\alpha_t
\big(
\nabla F_\infty(u_t)+r_t+\zeta_t
\big)
\right].
\]
Let $d_t:= \sup_{u\in\mathcal U}
\|\nabla F_t(u)-\nabla F_\infty(u)\|$. By Lemma~\ref{lem:grad_error},
\[
\|r_t\|
\le
\bar L e_t + d_t.
\]
Consequently, \eqref{eq:weighted_e_L1} and Lemma~\ref{lem:limiting_cost} give
\begin{equation}
\sum_{t=1}^{\infty}\alpha_t\mathbb E[\|r_t\|]\le \bar L\sum_{t=1}^{\infty}\alpha_t \mathbb{E}[e_t] + \sum_{t=1}^{\infty}\alpha_td_t<\infty.
\label{eq:weighted_r_L1}
\end{equation}

Let \(V_t:=\|u_t-u_\infty^\star\|^2\). By the optimality condition for \(u_\infty^\star\),
\[
u_\infty^\star
=
\Pi_{\mathcal U}
\left[
u_\infty^\star-\alpha_t\nabla F_\infty(u_\infty^\star)
\right],
\]
As a result, using the nonexpansiveness of the projection, we have
\begin{align}
V_{t+1}
&=
\|u_{t+1}-u_\infty^\star\|^2 \notag\\
&\le
\left\|
u_t-u_\infty^\star
-\alpha_t
\left(
\nabla F_\infty(u_t)-\nabla F_\infty(u_\infty^\star)
+r_t+\zeta_t
\right)
\right\|^2.
\label{eq:V_expand_start}
\end{align}
Expanding the right-hand side of \eqref{eq:V_expand_start} gives
\begin{align}
V_{t+1}
&\le
V_t
-2\alpha_t
\big(\nabla F_\infty(u_t)-\nabla F_\infty(u_\infty^\star)\big)^\top
(u_t-u_\infty^\star) \notag\\
&\quad
-2\alpha_t r_t^\top(u_t-u_\infty^\star)
-2\alpha_t \zeta_t^\top(u_t-u_\infty^\star)\notag\\
&\quad
+\alpha_t^2\|g_t-\nabla F_\infty(u_\infty^\star)\|^2.
\label{eq:V_expand}
\end{align}
Since \(F_\infty\) is \(\mu_\infty\)-strongly convex,
\begin{equation*}
\big(\nabla F_\infty(u_t)-\nabla F_\infty(u_\infty^\star)\big)^\top
(u_t-u_\infty^\star)
\ge
\mu_\infty V_t.
\end{equation*}
Moreover, \(\mathcal U\) is bounded, so there exists \(R>0\) such that $\|u_t-u_\infty^\star\|\le 2R$.
Therefore,
\begin{equation*}
-2\alpha_t r_t^\top(u_t-u_\infty^\star)
\le
4R\alpha_t\|r_t\|.
\end{equation*}
We now take conditional expectation with respect to \(\mathcal F_{t-1}\). By Lemma~\ref{lem:grad_error},
\begin{equation*}
\mathbb E[\zeta_t\mid\mathcal F_{t-1}]=0.
\end{equation*}
Since \(u_t-u_\infty^\star\) is \(\mathcal F_{t-1}\)-measurable, it follows that
\begin{equation*}
\mathbb E[
\zeta_t^\top(u_t-u_\infty^\star)
\mid
\mathcal F_{t-1}
]
=
\mathbb E[\zeta_t\mid\mathcal F_{t-1}]^\top
(u_t-u_\infty^\star)
=
0.
\end{equation*}
Finally, by the boundedness of the gradient estimate, there exists \(G_\infty>0\) such that
\[
\|g_t-\nabla F_\infty(u_\infty^\star)\|\le G_\infty.
\]
Taking conditional expectation in \eqref{eq:V_expand} yields
\begin{align}
\mathbb E[V_{t+1}\mid\mathcal F_{t-1}]
&\le
(1-2\mu_\infty\alpha_t)V_t\notag\\
&\quad+
4R\alpha_t\mathbb E[\|r_t\|\mid\mathcal F_{t-1}]
+
\alpha_t^2G_\infty^2.
\label{eq:RS_recursion}
\end{align}

Define
\[\varepsilon_t^{\rm RS}:=4R\alpha_t\mathbb E[\|r_t\|\mid\mathcal F_{t-1}]
+
\alpha_t^2G_\infty^2\]
and 
\[\vartheta_t^{\rm RS}:=2\mu_\infty\alpha_tV_t.\]
By the tower property, \eqref{eq:weighted_r_L1}, and Assumption~\ref{ass:stepsize},
\[
\sum_{t=1}^{\infty}\mathbb E[\varepsilon_t^{\rm RS}]=
4R\sum_{t=1}^{\infty}\alpha_t\mathbb E[\|r_t\|]
+
G_\infty^2\sum_{t=1}^{\infty}\alpha_t^2 <\infty.
\]
Since $\varepsilon_t^{\rm RS}\ge 0$, Tonelli's theorem implies $\mathbb E\left[\sum_{t=1}^{\infty}
\varepsilon_t^{\rm RS}\right] = \sum_{t=1}^{\infty}\mathbb E[\varepsilon_t^{\rm RS}]$. Hence,
\begin{equation}
\sum_{t=1}^{\infty}
\varepsilon_t^{\rm RS}
<\infty
\qquad
\text{almost surely}.
\label{eq:RS_error_sum}
\end{equation}

Set \(\mathcal G_t:=\mathcal F_{t-1}\). Then \(V_t\),
\(\varepsilon_t^{\rm RS}\), and \(\vartheta_t^{\rm RS}\) are
nonnegative and \(\mathcal G_t\)-measurable, and
\eqref{eq:RS_recursion} can be written as
\[
\mathbb E[V_{t+1}\mid\mathcal G_t]
\leq
V_t
-
\vartheta_t^{\rm RS}
+
\varepsilon_t^{\rm RS}.
\]
Thus, \eqref{eq:RS_error_sum} allows us to apply the
Robbins--Siegmund supermartingale convergence theorem
\cite[Lemma~4.1]{franci2022convergence}, with
\[
v_t=V_t,
\qquad
\delta_t=0,
\qquad
\varepsilon_t=\varepsilon_t^{\rm RS},
\qquad
\vartheta_t=\vartheta_t^{\rm RS}.
\]
It follows that \(V_t\) converges almost surely and
\[
\sum_{t=1}^{\infty}
\vartheta_t^{\rm RS}
=
2\mu_\infty
\sum_{t=1}^{\infty}\alpha_tV_t
<\infty
\qquad
\text{almost surely}.
\]
Since \(\sum_{t=1}^{\infty}\alpha_t=\infty\), we have $\liminf_{t\to\infty}V_t=0$.
Because \(V_t\) converges almost surely, its limit must be zero.
Therefore,
\[
u_t\to u_\infty^\star
\qquad
\text{almost surely}.
\]
\hfill\(\blacksquare\)

\bibliographystyle{IEEEtran}
\bibliography{ssfo.bib}

\end{document}